\documentclass{amsart}
\usepackage{amsmath,amssymb,amsthm,amscd}
\usepackage{graphicx}
\usepackage[all]{xy}
\theoremstyle{plain}
\newtheorem{thm}{Theorem}

\begin{document}

\title{On the arithmetic Chern character}
\author{H. \sc Gillet}
\email{gillet@uic.edu}
\address{Department of Mathematics, Statistics, and Computer Science\\
University of Illinois at Chicago\\
322 Science and Engineering Offices (M/C 249)\\
851 S. Morgan Street\\
Chicago, IL 60607-7045, USA}
\thanks{This material is based upon work supported in part by the National Science Foundation under Grant No. DMS-0901373}
\author{C. \sc Soul\'e}
\email{soule@ihes.fr}
\address{IH\'ES, 35 route de Chartres,\, 91440 Bures-sur-Yvette, France}
\date{\today }


\maketitle

Let $X$ be a proper and flat scheme over ${\mathbb Z}$, with smooth generic fiber $X_{\mathbb Q}$. In \cite{GS} we attached to every hermitian vector bundle $\overline E=(E,\|\;\|)$ on $X$ a Chern character class lying in the arithmetic Chow groups of $X$:
$$
\widehat{\rm ch} (\overline E) \in \underset{p \geq 0}{\oplus} \ \widehat{\rm CH}^p (X) \otimes {\mathbb Q} \, .
$$
Unlike the usual Chern character with values in the ordinary Chow groups, $\widehat{\rm ch}$ is not additive on exact sequences; indeed suppose that $\overline E_i$, $i=0,1,2$ is a triple of hermitian vector bundles on X, and that we are give an exact sequence
$$
0 \to E_0 \to E_1 \to E_2 \to 0
$$
of the underlying vector bundles on $X$, (\emph{i.e.} in which we ignore the hermitian metrics). Then the difference $\widehat{\rm ch} (\overline E_0) + \widehat{\rm ch} (\overline E_2) - \widehat{\rm ch} (\overline E_1)$, is represented by a secondary characteristic class $\widetilde{\rm ch}$ first introduced by Bott and Chern \cite{BC} and defined in general in \cite{Bi-G-S}.  These Bott-Chern forms measure the defect in additivity of the Chern forms associated by Chern-Weil theory to the hermitian bundles in the exact sequence.

\smallskip

Assume now that the sequence
$$
0 \to E_0 \to E_1 \to E_2 \to 0 \eqno (*)
$$
is exact on the generic fiber $X_{\mathbb Q}$ but not on the whole of $X$. We shall prove here (Theorem~1) that
$\widehat{\rm ch} (\overline E_0) + \widehat{\rm ch} (\overline E_2) - \widehat{\rm ch} (\overline E_1)$ is  the sum of the class of $\widetilde{\rm ch}$ and the localized Chern character of $(*)$ (see \cite{Fu}, 18.1). This result fits well with the idea that characteristic classes with support on the finite fibers of $X$ are the non-archimedean analogs of Bott-Chern classes (see \cite{GS3}).

\smallskip

In Theorem~2 we compute more explicitly these secondary characteristic classes in a situation encountered when proving a ``Kodaira vanishing theorem'' on arithmetic surfaces (\cite{S}, 3.3).

\bigskip

\noindent {\bf Notation.} If $A$ is an abelian group we let
$
A_{\mathbb Q} = A \underset{\mathbb Z}{\otimes} {\mathbb Q} \, .
$

\section{A general formula}\label{sec1}

\subsection{ \ } Let $S = {\rm Spec} ({\mathbb Z})$ and $f : X \to S$ a flat scheme of finite type over $S$. We assume that the generic fiber $X_{\mathbb Q}$ is smooth and equidimensional of dimension $d$. For every integer $p \geq 0$ we denote by $A^{pp} (X_{\mathbb R})$ the real vector space of smooth real differential forms $\alpha$ of type $(p,p)$ on the complex manifold $X({\mathbb C})$ such that $F_{\infty}^* (\alpha) = (-1)^p \, \alpha$, where $F_{\infty}$ is the  anti-holomorphic involution of $X({\mathbb C})$ induced by complex conjugation. Let
$$
A(X) = \underset{p \geq 0}{\oplus} \ A^{pp} (X_{\mathbb R})
$$
and
$$
\widetilde A (X) = \underset{p \geq 1}{\oplus} \ A^{p-1,p-1} (X_{\mathbb R}) / ({\rm Im} (\partial) + {\rm Im} (\overline\partial)) \, .
$$
For every $p \geq 0$ we let $\widehat{\rm CH}_p (X)$ be the $p$-th {\it arithmetic Chow homology} group of $X$ (\cite{ARR}, \S2.1, Definition 2). Elements  of $\widehat{\rm CH}_p (X)$ are represented by pairs $(Z,g)$ consisting of a $p$-dimensional cycle $Z$ on $X$ and a Green current $g$ for $Z({\mathbb C})$ on $X({\mathbb C})$. Recall that here a Green current for $Z({\mathbb C})$ is a current (\emph{i.e.} a form with distribution coefficients) of type $(p-1,p-1)$ such that $dd^c (g) + \delta_Z$ is $C^\infty$,  $\delta_Z$ being the current of integration on $Z({\mathbb C})$ ). There are canonical morphisms (\cite{ARR}, 2.2.1):
\begin{eqnarray}
z : \widehat{\rm CH}_p (X) &\to &{\rm CH}_p (X) \nonumber \\
(Z,g) &\mapsto &Z \nonumber
\end{eqnarray}
and
\begin{eqnarray}
\omega : \widehat{\rm CH}_p (X) &\to &A^{pp} (X_{\mathbb R}) \nonumber \\
(Z,g) &\mapsto &dd^c (g) + \delta_Z  \nonumber.
\end{eqnarray}

\smallskip

Let ${\rm CH}_p^{\rm fin} (X)$ be the Chow homology group of cycles on $X$ the support of which does not meet $X_{\mathbb Q}$. There is a canonical morphism
$$
b : {\rm CH}_p^{\rm fin} (X) \to \widehat{\rm CH}_p (X)
$$
mapping the class of $Z$ to the class of $(Z,0)$. The composite morphism
$$
z \circ b : {\rm CH}_p^{\rm fin} (X) \to {\rm CH}_p (X)
$$
is the obvious map. Let
$$
a : A^{d-p-1,d-p-1} (X_{\mathbb R}) \to \widehat{\rm CH}_p (X)
$$
be the map sending $\eta$ to the class of $(0,\eta)$. We have
$$
\omega \circ a (\eta) =  dd^c (\eta) \, .
$$

\subsection{ \ } We assume given a sequence
$$
0 \to \overline E_0 \to \overline E_1 \to \overline E_2 \to 0
$$
of hermitian vector bundles on $X$, the restriction of which to $X_{\mathbb Q}$ is exact. Let
$$
{\rm ch}^{\rm fin} (E_{\bullet}) \cap [X] \in {\rm CH}_{\rm fin} (X)_{\mathbb Q} = \underset{p \geq 0}{\oplus} \ {\rm CH}_p^{\rm fin} (X)_{\mathbb Q}
$$
be the {\it localized Chern character} of $E_{\bullet}$ (\cite{Fu} 18.1), and
$$
\widetilde{\rm ch} (\overline E_{\bullet}) \in \widetilde A (X)_{\mathbb Q}
$$
the Bott-Chern secondary characteristic class \cite{Bi-G-S}, such that
$$
dd^c \, \widetilde{\rm ch} (\overline E_{\bullet}) = \sum_{i=0}^2 (-1)^i \, {\rm ch} (\overline E_{i,{\mathbb C}}) \, ,
$$
where ${\rm ch} (\overline E_{i,{\mathbb C}}) \in A(X)$ is the differential form representing the Chern character of the restriction $E_{i,{\mathbb C}}$ of $E_i$ to $X({\mathbb C})$. Finally, if $i=0,1,2$, we let
$$
\widehat{\rm ch} (\overline E_i) \cap [X] \in \widehat{\rm CH} (X)_{\mathbb Q} = \underset{p \geq 0}{\oplus} \ \widehat{\rm CH}_p (X)_{\mathbb Q}
$$
be the {\it arithmetic Chern character} of $\overline E_i$ (\cite{GS} 4.1,\, \cite{ARR} Theorem 4).
\begin{thm}\label{thm1}
The following equality holds in $\widehat{\rm CH} (X)_{\mathbb Q}$:
$$
\sum_{i=0}^2 (-1)^i \, \widehat{\rm ch} (\overline E_i) \cap [X] = b \, ({\rm ch}^{\rm fin} (E_{\bullet}) \cap [X]) + a \, (\widetilde{\rm ch} (\overline E_{\bullet})) \, .
$$
\end{thm}
\subsection{ \ }  
This theorem is a special case of Lemma 21 in \cite{ARR}, though this may not be immediately apparent. Therefore, for the sake of completeness, we give a proof here.
\subsection{ \ } To prove Theorem~\ref{thm1} we consider the Grassmannian graph construction applied to $E_{\bullet}$ (\cite{Fu} 18.1, \cite{ARR} 1.1). It consists of a proper surjective map
$$
\pi : W \to X \times {\mathbb P}^1
$$
such that, if $\phi \subset X$ is the support of the homology of $E_{\bullet}$ (hence $\phi_{\mathbb Q}$ is empty), the restriction of $\pi$ onto $(X -\phi) \times {\mathbb P}^1$ and $X \times {\mathbb A}^1$ is an isomorphism. The effective Cartier divisor
$$
W_{\infty} = \pi^{-1} (X \times \{ \infty \} )
$$
is the union of the Zariski closure $\widetilde X$ of $(X - \phi) \times \{ \infty \}$ with $Y = \pi^{-1} (\phi \times \{ \infty \})$. The sequence $E_{\bullet}$ extends to a complex
$$
0 \to \widetilde E_0 \to \widetilde E_1 \to \widetilde E_2 \to 0 \, ,
$$
which is isomorphic to the pull-back of $E_{\bullet}$ over $X \times {\mathbb A}^1$. The restriction of $\widetilde E_{\bullet}$ to $\widetilde X$ is canonically split exact. On $W_{\mathbb Q} = X_{\mathbb Q} \times {\mathbb P}_{\mathbb Q}^1$ the sequence $\widetilde E_{\bullet}$ is exact; it coincides with $E_{\bullet}$ (resp. $0 \to E_0 \to E_0 \oplus E_2 \to E_2 \to 0$) when restricted to $X_{\mathbb Q} \times \{ 0 \}$ (resp. $X_{\mathbb Q} \times \{ \infty \}$). We choose a metric on $\widetilde E_{\bullet}$ for which these isomorphisms are isometries.

\subsection{ \ } Let
$$
x = \sum_{i=0}^2 (-1)^i \, \widehat{\rm ch} (\overline{\widetilde E}_i) \, ,
$$
and denote by $t$ the standard parameter of ${\mathbb A}^1$. In the arithmetic Chow homology of $W$ we have
$$
0 = x \cap (W_0 - W_{\infty} , - \log \vert t \vert^2) \, .
$$
If $x$ is the class of $(Z,g)$, with $Z$ meeting properly $W_0$ and $W_{\infty}$, we get
$$
x \cap (W_0 - W_{\infty} , -\log \vert t \vert^2) = (Z \cap (W_0 - W_{\infty}) , g * (-\log \vert t \vert^2) )\, ,
$$
where the $*$-product is equal to
$$
g * (-\log \vert t \vert^2) = g (\delta_{W_0} - \delta_{W_{\infty}}) - {\rm ch} (\overline{\widetilde E}_{\bullet}) \log \vert t \vert^2 \, .
$$
Since $W_{\infty} = \widetilde X \cup Y$, with $Y_{\mathbb Q} = \emptyset$, we get
\begin{eqnarray}
\label{eq1}
0 &= &x \cap (W_0 - W_{\infty} , - \log \vert t \vert^2) \\
&= &(Z \cap W_0 , g \, \delta_{W_0}) - (Z \cap \widetilde X , g \, \delta_{\widetilde X}) - (Z \cap Y , 0) - (0,{\rm ch} (\overline{\widetilde E}_{\bullet}) \log \vert t \vert^2) \, .  \nonumber
\end{eqnarray}
The restriction of $\overline{\widetilde E}_{\bullet}$ to $\widetilde X$ is split exact, therefore
$$
(Z \cap \widetilde X , g \, \delta_{\widetilde X}) = 0 \, .
$$
Applying $\pi_*$ to (\ref{eq1}) we get
\begin{equation}
\label{eq2}
0 = \widehat{\rm ch} (\overline E_{\bullet}) - \pi_* (Z \cap Y) , 0) - (0,\pi_* ({\rm ch} (\overline{\widetilde E}_{\bullet}) \log \vert t \vert^2) \, .
\end{equation}
By definition of the localized Chern character (\cite{Fu}, 18.1, (14))
\begin{equation}
\label{eq3}
\pi_* (Z \cap Y) = {\rm ch}^{\rm fin} (E_{\bullet}) \cap [X]
\end{equation}
in ${\rm CH}^{\rm fin} (X)_{\mathbb Q}$. On the other hand we deduce from \cite{GS}, (1.2.3.1), (1.2.3.2) that
\begin{equation}
\label{eq3b}
- \, \pi_* ({\rm ch} (\overline{\widetilde E}_1) \log \vert t \vert^2) = \widetilde{\rm ch} (\overline{E}_{\bullet}) \, .
\end{equation}
and upon replacing $t$ by $1/t$, as in the proof of (1.3.2) in \cite{GS}, we see that
\begin{equation}
\label{eq4}
\pi_* ({\rm ch} (\overline{\widetilde E}_{\bullet}) \log \vert t \vert^2) = - \, \pi_* ({\rm ch} (\overline{\widetilde E}_1) \log \vert t \vert^2)\;.
\end{equation}
Theorem~\ref{thm1} follows from (\ref{eq2}), (\ref{eq3}), (\ref{eq3b}), (\ref{eq4}).

\section{A special case}\label{sec2}

\subsection{ \ } We keep the hypotheses of the previous section, and we assume that $X$ is normal,
$d=1$,  $E_0$ and $E_2$ have rank one and the metrics on $E_0$ and $E_2$ are induced by the metric on $E_1$. Finally, we assume that there exists a closed subscheme $\phi$ in $X$ which is $0$-dimensional and such that there is an exact sequence of sheaves on $X$
\begin{equation}
\label{eq5}
0 \to E_0 \to E_1 \to E_2 \otimes I_{\phi} \to 0 \, ,
\end{equation}
where $I_{\phi}$ is the ideal of definition of $\phi$.

\smallskip

Let
$$
\widetilde c_2 \in A^{1,1} (X_{\mathbb R}) / ({\rm Im} (\partial) + {\rm Im} (\overline\partial))
$$
be the second Bott-Chern class of (\ref{eq5}), $\Gamma (\phi , {\mathcal O}_{\phi})$ the finite ring of functions on $\phi$ and $\# \, \Gamma (\phi , {\mathcal O}_{\phi})$ its order. Let
$$
f_* : \widehat{\rm CH}_0 (X)_{\mathbb Q} \to \widehat{\rm CH}_0 (S) = {\mathbb R}
$$
be the direct image morphism.

\begin{thm}\label{thm2}
We have an equality of real numbers
$$
f_* (\widehat c_2 (\overline E_1) \cap [X]) = f_* (\widehat c_1 (\overline E_0) \, \widehat c_1 (\overline E_2) \cap [X]) - \int_{X({\mathbb C})} \widetilde c_2 + \log \# \, \Gamma (\phi , {\mathcal O}_{\phi}) \, .
$$
\end{thm}

\subsection{ \ } To prove Theorem~\ref{thm2} we remark first that
$$
\widehat c_1 (\overline E_1) = \widehat c_1 (\overline E_0) + \widehat c_1 (\overline E_2) \, ,
$$
because the metrics on $E_0$ and $E_2$ are induced from $\overline E_1$. Therefore, since ${\rm ch}_2 = -c_2 + \frac{c_1^2}{2}$, we get
\begin{eqnarray}
\widehat{\rm ch}_2 (\overline E_1) &= & - \widehat c_2 (\overline E_1) + \frac{(\widehat c_1 (\overline E_0) + \widehat c_1 (\overline E_2))^2}{2} \nonumber \\
&= &- \widehat c_2 (\overline E_1) + c_1 (\overline E_0) \, \widehat c_1 (\overline E_2) + \widehat{\rm ch_2} (\overline E_0) + \widehat{\rm ch_2} (\overline E_2) \, . \nonumber
\end{eqnarray}
By Theorem~\ref{thm1}, this implies that
\begin{equation}
\label{eq6}
\widehat c_2 (\overline E_1) \cap [X] = \widehat c_1 (\overline E_0) \, \widehat c_1 (\overline E_2) \cap [X] + b \, ({\rm ch}^{\rm fin} (E_{\bullet}) \cap [X]) + a \, (\widetilde{\rm ch} (\overline E_{\bullet})) \, .
\end{equation}
Since $\widetilde{\rm ch}_0 (\overline E_{\bullet})$ and $\widetilde{\rm ch}_1 (\overline E_{\bullet})$ vanish we have
$$
\widetilde{\rm ch} (\overline E_{\bullet}) = - \, \widetilde c_2 \, .
$$

Therefore, if we apply $f_*$ to (\ref{eq6}), we get
$$
f_* (\widehat c_2 (\overline E_1) \cap [X]) = f_* (\widehat c_1 (\overline E_0) \, \widehat c_1 (\overline E_2) \cap [X]) - \int_{X({\mathbb C})} \widetilde c_2 + f_* (b ({\rm ch}^{\rm fin} (E_{\bullet}) \cap [X])) \, ,
$$
and we are left with showing that
\begin{equation}
\label{eq7}
f_* \circ b \, ({\rm ch}^{\rm fin} (E_{\bullet}) \cap [X]) = \log \# \, \Gamma (\phi , {\mathcal O}_{\phi}) \, .
\end{equation}
Let $\vert \phi \vert = \{ P_1,\cdots ,P_n \} \subset X$ be the support of $\phi$ and $\psi = f(\vert \phi \vert ) \subset S$. The following diagram is commutative:
$$
\xymatrix{
{\rm CH}_0 (\phi) \ar[d]^{f_*}\ar[r]^b &\widehat{\rm CH}_0 (X) \ar[d]^{f_*} \\
{\rm CH}_0 (\psi) \ar[r]^{b \quad } &\widehat{\rm CH}_0 (S) = {\mathbb R} \, , \\
}
$$
where
$$
b : {\rm CH}_0 (\psi) = {\mathbb Z}^{\psi} \to {\mathbb R}
$$
maps $(n_p)_{p \in \psi}$ to $\underset{p}{\sum} \, n_p \log (p)$.

\smallskip

For any prime $p \in \psi$ we let ${\mathbb Z}_{(p)}$ be the local ring of $S$ at $p$ and we let $\ell_p = \ell_p (\phi)\geq 0$ be the length of the finite ${\mathbb Z}_{(p)}$-module $\Gamma (\phi , {\mathcal O}_{\phi}) \otimes {\mathbb Z}_{(p)}$. Clearly
$$
\log \# \, \Gamma (\phi , {\mathcal O}_{\phi}) = \sum_{p \in \psi} \ell_p \log (p) \, ,
$$
hence it is enough to prove that
\begin{equation}
\label{eq8}
f_* ({\rm ch}^{\rm fin} (E_{\bullet}) \cap [X]) = (\ell_p) \in {\rm CH}_0 (\psi)_{\mathbb Q} = {\mathbb Q}^{\psi} \, .
\end{equation}

The complex $E_{\bullet}$ defines an element
$$
[E_{\bullet}] = \sum_{i=1}^{n} \, [{\mathcal O}_{\phi , P_i}] \in K_0^{\phi} (X) = \bigoplus_{i=1}^{n} K_0^{P_i} (X) \, .
$$
To prove (9), by replacing $X$ by an affine neighbourhood of $P$, one can assume that $\vert \phi \vert = \{ P \}$, and it is enough to show that, if
$
p = f(P) \, ,
$

$$
f_* ({\rm ch}^{\rm fin} ({\mathcal O} _{\phi , P}) \cap [X]) = \ell_p ({\mathcal O}_{\phi , P}) [p] \, .
$$

Now recall that, if ${\mathcal F}$ is a coherent sheaf on a scheme $X$ of finite type over $S$, supported on a finite set of closed points, the associated $0$-cycle
$$
[{\mathcal F}] = \sum_{P \in \vert {\mathcal F} \vert} \ell_p ({\mathcal F}_P)[P] \in Z_0 (X)
$$
is such that, if $f : X \to Y$ is a proper morphism of schemes of finite type over $S$,
$$
f_* [{\mathcal F}] = [f_* ({\mathcal F})]
$$
(\cite{Fu} , 15.1.5). Hence it is enough to show that
\begin{equation}
\label{eq10}
{\rm ch}^P ({\mathcal O}_{\phi , P}) = \ell_p ({\mathcal O}_{\phi , P}) [P] \in {\rm CH}_0 (P)_{\mathbb Q} \simeq {\mathbb Q} \, .
\end{equation}
Replacing $X$ by an affine neighbourhood of $P$, we may assume that we have an exact sequence
\begin{equation}
\label{eq11}
0 \longrightarrow {\mathcal O}_X \overset{\alpha}{\longrightarrow} {\mathcal O}_X^2 \overset{\beta}{\longrightarrow} {\mathcal O}_X 
\longrightarrow {\mathcal O}_{\phi} \longrightarrow 0 \, .
\end{equation}

Hence the ideal $I_{\phi} \subset {\mathcal O}_X (X)$ is generated by two elements $\beta_1$ and $\beta_2$. Since $X$ is normal, its local rings satisfy Serre property $S_2$ and, as $\dim (X) = 2$, $X$ is Cohen-Macaulay. Since $\beta_1$ and $\beta_2$ span an ideal of height two, $(\beta_1 , \beta_2)$ is a regular sequence and the sequence (\ref{eq11}) is isomorphic to the Koszul resolution of ${\mathcal O}_{\phi} = {\mathcal O}_X / (\beta_1 , \beta_2)$. Now (\ref{eq10}) can be deduced from the following general fact:

\bigskip

\noindent {\bf Lemma 1.} {\it Let $X = {\rm Spec} \, (A)$ be an affine scheme and $Z \subset X$ a closed subset such that the ideal $I_Z = (x_1 , \ldots , x_n)$ is generated by a regular sequence $(x_1 , \ldots , x_n)$. Let $K_{\bullet} (x_1 , \ldots , x_n)$ be the Koszul complex associated to $(x_1 , \ldots , x_n)$. Then
$$
{\rm ch}_n^Z (K_{\bullet} (x_1 , \ldots , x_n)) = [{\mathcal O}_Z] \in {\rm CH}_0 (Z)_{\mathbb Q} \, .
$$
}

\bigskip

\noindent {\bf Proof.} The Grassmannian-graph construction on $K_{\bullet} (x_1 , \ldots , x_n)$ coincides with the deformation to the normal bundle of $Z$ in $X$. If $W$ is defined as in 1.4,
$$
W_{\infty} = \widetilde X \cup \widehat{\mathbb P} (N_{Z/X}) \, ,
$$
where $\widetilde X$ is the blow up of $X$ along $Z$, and $\widehat P (N_{Z/X})$ is the projective completion of the normal bundle of $Z$ in $X$. The pull back of the Koszul complex $K_{\bullet} (x_1 , \ldots , x_n)$ to $W \backslash W_{\infty}$ extends to a complex $\widetilde K_{\bullet} (x_1 , \ldots , x_n)$ on $W$. The restriction of $\widetilde K_{\bullet} (x_1 , \ldots , x_n)$ to $\widetilde X$ is acyclic while the restriction of $\widetilde K_{\bullet} (x_1 , \ldots  , x_n)$ to $\widehat{\mathbb P} (N_{Z/X})$ is a resolution of the structure sheaf of the zero section $Z \subset N_{Z/X} \subset \widehat{\mathbb P} (N_{Z/X})$.

\smallskip

Now observe that $Z \subset \widehat{\mathbb P} (N_{Z/X})$ is an intersection of Cartier divisors $D_1 , \ldots ,$ $D_n$, hence
\begin{eqnarray}
&&{\rm ch} (\widetilde K_{\bullet} (x_1 , \ldots , x_n) \mid_{\widehat{\mathbb P} (N_{Z/X})}) \nonumber \\
&= &\prod_{i=1}^{n} {\rm ch} ({\mathcal O}(-D_i) \to {\mathcal O}_{\widehat{\mathbb P} (N_{Z/X})}) \nonumber \\
&= &\prod_{i=1}^{n} {\rm ch} ({\mathcal O}(D_i)) \, . \nonumber
\end{eqnarray}
Since
$$
{\rm ch} ({\mathcal O}_{D_i}) = {\rm ch}_1 ({\mathcal O}_{D_i}) + x_i = [D_i] + x_i
$$
where $x_i$ has degree $\geq 2$, we get
$$
{\rm ch} (\widetilde K_{\bullet} (x_1 , \ldots , x_n) \mid_{\widehat{\mathbb P} (N_{Z/X})}) = [D_1] \ldots [D_n] = [Z] \, .
$$
This ends the proof of Lemma~1 and Theorem~2.

\vglue 1cm

\end{document}